\documentclass[12pt]{amsart}
\newtheorem{theorem}{Theorem}
\newtheorem{lemma}[theorem]{Lemma}
\newtheorem{corollary}[theorem]{Corollary}

\newcommand{\N}{{\Bbb N}}

\newcounter{rmnum}

\newcommand{\cone}{\operatorname{Cone}}
\newcommand{\e}{\operatorname{e-dim}}

\newenvironment{alphanum}{\begin{list}{{\rm (\alph{rmnum})}}
{\usecounter{rmnum}\def\makelabel##1{\hss\llap{##1}}
\setlength{\leftmargin}{37pt}\setlength{\itemindent}{0pt}
\setlength{\topsep}{5pt}\setlength{\parsep}{0pt}\setlength
{\itemsep}{0pt}}}{\end{list}}

\begin{document}
\title{Universal Metric Spaces and Extension Dimension}
\author{Alex Chigogidze and Vesko Valov\\[3pt] }
\address{Department of Mathematics and Statistics, 
University of Saskatchewan, McLean Hall,
106 Wiggins Road, Saskatoon, SK, Canada S7N 5E6}
\email{chigogid@math.usask.ca}

\bigskip
\address{Department of Mathematics,
University of Swaziland, Pr. Bag 4,
Kwaluseni, Swaziland}
\email{valov@realnet.co.sz}


\begin{abstract}
For any countable $CW$-complex $K$ and a cardinal number
$\tau\geq\omega$ we construct a completely
metrizable space $X(K,\tau)$ of weight $\tau$ with the
following properties:
$\e X(K,\tau)\leq K$, $X(K,\tau)$ is an absolute extensor
for all normal spaces
$Y$ with $\e Y\leq K$, and for any completely metrizable
space $Z$ of weight $\leq\tau$ 
and $\e Z\leq K$ the set of closed embeddings
$Z\rightarrow X(K,\tau)$ is dense in
the space $C(Z,X(K,\tau))$ of all continuous maps from
$Z$ into $X(K,\tau)$ endowed
with the limitation topology. This result is applied to prove
the existence of universal spaces for all metrizable spaces of
given weight and with a given cohomological dimension.
\end{abstract}

\maketitle

\bigskip
The existence of universal separable metric spaces for
extension dimension with respect to countable $CW$-complexes was
proved by Olszewski in \cite{wo:95}. In the class of all 
metric spaces of a given weight this problem was recently solved
by Levin \cite{l:98}. In the present note we show the existence
of universal metric spaces having some extra properties (see Theorem 1 below).  
The concept of extension dimension was introduced by
Dranishnikov \cite{d:94}
(see also \cite{ch1:97}, \cite{dd:96}). For a normal space
$X$ and a
$CW$-complex $K$ we write $\e X\leq K$ (the extension
dimension of $X$
does not exceed $K$) if $K$ is an absolute extensor for
$X$. This means that any continuous map $f \colon A \to K$,
defined on a closed subset $A$ of $X$, admits a continuous
extension $\bar{f} \colon X \to K$.
Since not every $CW$-complex is an absolute neighborhood extensor
for normal spaces, we can enlarge the class of normal spaces $X$
with $\e X\leq K$ ($K$ is a $CW$-complex) by introducing the following
notion (see \cite[Definition 2.5]{lr:91}): A normal space $X$ is in the 
class $\alpha (K)$ if every continuous map from a closed 
$A\subset X$ to $K$ which extends to a map of a neighborhood of $A$
to $K$ can be extended to a map of $X$ to $K$. Obviously, if
$K\in ANE(X)$ (this, for example, holds for every $X$ admitting a
perfect map onto a first countable paracompact space \cite{jd:95})
then $X\in\alpha (K)$ if and only if $\e X\leq K$.
We also adopt the following definition: a continuous map
$f: X\rightarrow Y$ is called
$K$-soft (resp., $K$-soft with respect to metrizable spaces) if for
any normal (resp., metrizable) space $Z$ with $Z\in\alpha (K)$,
any closed $Z_0\subset Z$, and 
any two maps $g: Z_0\rightarrow X$, $h: Z\rightarrow Y$ with
$f\circ g=h|Z_0$, there exists a map 
$k: Z\rightarrow X$ such that $f\circ k=h$ and $k|Z_0=g$.  

For any $CW$-complex $K$ and any cardinal number $\tau\geq\omega$
let ${\mathcal M}(K,\tau)$ be the 
class of all completely metrizable spaces $X$ of weight
$\tau$ with $\e X\leq K$. 
The following theorem is our main result:

\begin{theorem}
For any countable $CW$-complex $K$ and a cardinal number
$\tau\geq\omega$ there exists a completely
metrizable space $X(K,\tau)$ and a $K$-soft map
$f(K,\tau): X(K,\tau)\rightarrow l_2(\tau)$
satisfying the following properties:

\begin{alphanum}
\item $X(K,\tau)\in {\mathcal M}(K,\tau)$.

\item $X(K,\tau)$ is an absolute extensor for all normal
spaces $Y$ with $Y\in\alpha (K)$.

\item $f(K,\tau)$ is strongly $(K,\tau)$-universal, i.e.
for any open cover $\mathcal U$ of
$X(K,\tau)$, any (complete) metric space $Z$ of weight $\leq\tau$ 
with $\e Z\leq K$ and any
map $g: Z\rightarrow X(K,\tau)$
there exists a (closed) embedding $h: Z\rightarrow X(K,\tau)$
$\mathcal U$-close to $g$ with
$f(K,\tau)\circ g=f(K,\tau)\circ h$.   
\end{alphanum}
\end{theorem}

For the case $K=S^n$, $n\in\N$, Theorem 1 was proved in
\cite[Theorem 2.7]{cv:90}.
Our proof of Theorem 1 is based on the next few lemmas and the
techniques developed in  
\cite{ch2:89} and \cite{cv:90}.

\begin{lemma}\label{L:1}
For any countable $CW$-complex $K$ and any
separable (completely) \linebreak
metrizable space $X$ there exists a separable (completely)
metrizable
space $Y_X$ with $\e Y_X\leq K$ and a $K$-soft map
$f: Y_X\rightarrow X$. 
\end{lemma}

\begin{proof}
Let $P$ be a Polish $ANR$ homotopically equivalent to $K$ and
$\varphi :X\rightarrow P$ and $\psi :P\rightarrow X$ be two maps
such that $\psi\circ\varphi$ is homotopic to $id_P$ and
$\psi\varphi$ is homotopic to $id_K$. For
extension dimension with respect to $P$ this lemma was proved in
\cite[Proposition 5.9]{ch1:97}. So, for a given (complete) separable  
metric space $X$ there is a (complete) separable metric space $Y_X$
with $\e Y_X\leq P$ and a $P$-soft map $f: Y_X\rightarrow X$.
According to next claim, $f$ is $K$-soft.  

\medskip\noindent
{\it Claim.} {\em If $Z\in\alpha (K)$ is normal, then $Z\in\alpha (P)$}

Suppose $Z\in\alpha (K)$
is a normal space. Since every Polish $ANR$ is an $ANE$ for normal spaces,
$Z\in\alpha (P)$ is equivalent to $P\in\AE(Z)$. Take a map
$g: A\rightarrow P$, where $A\subset Z$ is closed and consider the map
$\psi\circ g: A\rightarrow K$. Because $g$ can be extended to a map from a
neighborhood $U$ of $A$ into $P$, $\psi\circ g$ can be extended to a map from
$U$ to $K$. Since $Z\in\alpha (K)$, there is an extension $h: Z\rightarrow K$
of $\psi\circ g$. Then the restriction $(\varphi h)|A$ is homotopic to $g$.
Finally, using that the Homotopy Extension Theorem holds for normal spaces
and Polish $ANR$'s, we conclude that $g$ is extendable to a map from $Z$ into
$P$. Hence $Z\in\alpha (P)$.

\medskip
It remains only to show that $\e Y_X\leq K$. And this follows from 
$\e Y_X\leq P$ and the fact that the Homotopy Extension Theorem holds for
metric spaces and $CW$-complexes \cite{jd:95}.
\end{proof}

\begin{lemma}\label{L:2}
\cite{l:98} Let $f: X\rightarrow Y$ be a uniformly
0-dimensional map of metrizable spaces $X$ and $Y$.
Then $\e X\leq \e Y$.
\end{lemma}

Recall that a map $f: X\rightarrow Y$, where $X$ and $Y$ are
metrizable spaces, is called 
uniformly 0-dimensional \cite{k:52} if there exists a compatible
metric on $X$ such that 
for every $\varepsilon >0$ and every $y\in f(X)$ there is an open
neighborhood $U$ of $y$ such
that $f^{-1}(U)$ can be represented as the union of disjoint
open sets of $diam <\varepsilon$.
It is well known that every metric space admits a uniformly
0-dimensional map into Hilbert 
cube $Q$.  

\begin{lemma}\label{L:3}
For any countable $CW$-complex $K$ and a (completely) metrizable
space $Y$ of weight $\tau$ 
there exist a (completely) metrizable space X of weight
$\tau$ and a $K$-soft map 
$f: X\rightarrow Y$ such that $\e X\leq K$.
\end{lemma}

\begin{proof}
It suffices to prove this corollary when $Y$ is the space
$l_2(\tau)$. Fix a compatible 
metric $d_1$ on $l_2(\tau)$ and an uniformly 0-dimensional
surjection (with respect to 
$d_1$)
$g: l_2(\tau)\rightarrow A$ with $A$ a separable metric space.
By Lemma \ref{L:1}, there exists
a separable metric space $Z$ with $\e Z\leq K$ and a
$K$-soft map $h: Z\rightarrow A$.
Let $X$ be the fibered product of $l_2(\tau)$ and $Z$ with
respect to $g$ and $h$, and
let $f: X\rightarrow l_2(\tau)$ and $p: X\rightarrow Z$ denote
the corresponding 
projections of this fibered product. If $d_2$ is any metric on
$Z$, then $p$ is 
uniformly 0-dimensional with respect to the metric
$\displaystyle (d_1^2+d_2^2)^{1/2}$ 
(see \cite{ap:73}). Hence, by Lemma \ref{L:2}, $\e Z\leq K$. The
$K$-softness of $h$ implies
that $f$ is also $K$-soft. 

It remains only to show that $X$ is completely metrizable. To
this end let $B_X$ be the space obtained from $\beta X$ by
making the points of 
$\beta X-X$ isolated. According to \cite[Lemma 2]{p:78},
$B_X$ is paracompact,
and obviously, $B_X$ is first countable.

\medskip\noindent
{\it Claim.} {\em $\e B_X\leq K$}

Let $s: F\rightarrow K$ be an
arbitrary map, where $F\subset B_X$ is closed. Since
$\e X\leq K$, there exists 
an extension $s_1 \colon F\cup X\rightarrow K$ of $s$. Now we
need the following result
\cite[Theorem 11.2]{jd:95}: any contractible $CW$-complex is an
absolute
extensor for all
spaces admitting a perfect map onto a first countable paracompact
space.
This statement implies that $\cone(K)$ is an absolute extensor
for $B_X$.
Therefore, there exists an extension
$s_2: B_X\rightarrow \cone(K)$ of $s_1$. Let
$H=s_2^{-1}(\cone(K)-\{ b\} )$, where
$b\in \cone(K)-K$. Fix a retraction $r \colon \cone(K)-\{ b\}\rightarrow K$.
Since $H$ is clopen in $B_X$, it follows that $r\circ s_2$ can be
extended to a map $s_{3} \colon B_{X} \rightarrow K$.   
Then $s_3$ is an extension of $s$ and, consequently, $\e B_{X}\leq K$.

\medskip
Now let us go back to the proof of the completeness of $X$. Considering $X$
as a closed subset
of $B_X$ and using the fact that $l_2(\tau)$ is an absolute extensor for
paracompact spaces, we
can find a map $q \colon B_{X} \rightarrow l_2(\tau)$ such that $q|X=f$. Then,
since $f$ is $K$-soft
and $\e B_X\leq K$, there exists a retraction from $B_X$ onto $X$.
Finally, applying the
argument from \cite[the proof of Lemma 2]{p:78}, we conclude that $X$
is complete.
\end{proof}     

\medskip\noindent
{\it Proof of Theorem 1.} We will construct
an inverse sequence $S=\{X_n, p^{n+1}_{n}, n\in\N\}$ such that:

\begin{enumerate}
\item $X_{1}=l_2(\tau)$ and $X_n\in {\mathcal M}(K,\tau)$ for each
$n\geq 1$;

\item each $p^{n+1}_{n}: X_{n+1}\rightarrow X_n$ is a $K$-soft map
such that
for any completely metrizable space $Z$ of weight $\leq\tau$ with
$\e Z\leq K$
and any map $g: Z\rightarrow X_n$ there exists a closed embedding
$h: Z\rightarrow X_{n+1}$ with $p^{n+1}_{n}\circ h=g$.
\end{enumerate}

If $X_i$ and $p^{i}_{i-1}$ have already been constructed for $i=1,2,\dots , n$,
then, by Lemma \ref{L:3}, there exist a completely metrizable space $X_{n+1}$ of
weight $\tau$ and a $K$-soft map
$h_{n+1}\colon X_{n+1}\rightarrow X_n\times l_2(\tau)$
such that $\e X_{n+1}\leq K$. Let $p^{n+1}_{n}=\pi _n\circ h_{n+1}$, where
$\pi_{n} \colon X_{n} \times l_{2}(\tau)\rightarrow X_{n}$ is the natural projection.
Denote by $X(K,\tau)$ the limit space of $S$ and by $f(K,\tau)$ the
limit projection $p_{1}: X(K,\tau)\rightarrow X_{1}$. Obviously, $X(K,\tau)$
is a completely metrizable space of weight $\tau$ and $f(K,\tau)$ is $K$-soft. 
Following the proof of Lemma 2.6 from \cite{cv:90} one can show that
$f(K,\tau)$ is strongly $(K,\tau)$-universal.
Finally, by the limit theorem of Rubin-Schapiro \cite{rs:98},
$\e X(K,\tau)\leq K$. Observe that $X(K,\tau )$ is an absolute extensor for
all normal spaces $Y$ with $\e Y\leq K$ because $l_2(\tau)$ is an absolute
extensor for normal spaces and $f(K,\tau)$ is $K$-soft. 

We can apply Theorem 1 to obtain universal spaces for all metrizable spaces
with a given cohomological dimension and a given weight. Recall that, for
any abelian group $G$ and a natural number $n$, the cohomological dimension
$dim_GX$ of $X$ with a coefficient group $G$ is $\leq n$ iff
$\e X\leq K(G,n)$, where $X$ is a normal space and $K(G,n)$ is the
Eilenberg-MacLane complex. Let us agree the following notations: a map $f$
is called $(G,n)$-soft iff it is $K(G,n)$-soft, and $f$ is strongly
$(G,n,\tau)$-universal iff $f$ is strongly $(K(G,n),\tau)$-universal.

\begin{corollary}\label{Cor:1}
Let $G$ be a countable (resp., torsion) abelian group. Then for every
$n\in\N$ and $\tau\geq\omega$ there exists a completely metrizable space
$X_{\tau}(G,n)$ of weight $\tau$ and a map
$f_{\tau}(G,n): X_{\tau}(G,n)\rightarrow l_2(\tau)$ such that:

\begin{alphanum}
\item $\operatorname{dim}_GX_{\tau}(G,n)=n$.

\item $X_{\tau}(G,n)$ is an absolute extensor for all normal
(rep., metrizable) spaces $Y$ with $\operatorname{dim}_G Y\leq n$.

\item $f_{\tau}(G,n)$ is strongly $(G,n,\tau)$-universal and $(G,n)$-soft
(resp., $(G,n)$-soft with respect to metrizable spaces).
\end{alphanum}
\end{corollary}

\begin{proof}
If $G$ is countable, the proof follows directly from Theorem 1 with $K=K(G,n)$.
If $G$ is torsion, by \cite[Theorem B(a)]{jd:93}, there exists a countable
family $\sigma (G)$ of countable groups such that
$\operatorname{dim}_GY = \max\{\operatorname{dim}_HY:H\in\sigma (G)\}$ for any metrizable space $Y$. Then,
according to \cite[Lemma 2.4]{jd:96}, for each $n\in\N$ there is a countable
complex $K_n$ with $\operatorname{dim}_GY\leq n$ if and only if $\e Y\leq K_n$, $Y$ is any
metrizable space. Finally, apply Theorem 1 to $K_n$.
\end{proof}

Similarly, using Theorem 1 and \cite[Theorem B(a),(b) and (f)]{jd:93}), we
can obtain the following

\begin{corollary}\label{Cor:2}
Let $G$ be an arbitrary abelian group. Then for every $n\in\N$ and
$\tau\geq\omega$ there exists a completely metrizable space
$Y_{\tau}(G,n)$ of weight $\tau$ and a map
$g_{\tau}(G,n): Y_{\tau}(G,n)\rightarrow l_2(\tau)$ such that:

\begin{alphanum}
\item $\operatorname{dim}_GY_{\tau}(G,n)\leq n+1$.

\item $Y_{\tau}(G,n)$ is an absolute extensor for all metrizable
spaces $Z$ with $\operatorname{dim}_G Z\leq n$.

\item $g_{\tau}(G,n)$ is strongly $(G,n,\tau)$-universal and $(G,n)$-soft
with respect to metrizable spaces.
\end{alphanum}
\end{corollary}
    
For $\tau =\omega$ weaker versions of Corollary \ref{Cor:1} and
Corollary \ref{Cor:2} were proved in \cite{wo:95}.

\bibliographystyle{amsplain}

\bibliographystyle{amsplain}
\bibliography{triquot}

\end{document}